\documentclass[12pt]{article}
\input{amssym.def}
\input{amssym.tex}

\newcommand{\N}{{\Bbb N}}

\newcommand{\Z}{{\Bbb Z}}

\newcommand{\ga}{\gamma}

\newcommand{\pfend}{$/\!\! /$\bigskip}

\newtheorem{defi}{Definition}[section]
\newtheorem{theorem}{Theorem}

\newtheorem{lemma}[defi]{Lemma}
\newtheorem{proposition}[defi]{Proposition}

\newtheorem{corollary}[defi]{Corollary}

\newtheorem{remark}[defi]{Remark}

\newtheorem{*example}[defi]{*Example}

\addtocontents{ences.txt}{References}

\begin{document}

\title{\bf A Finitely Presented Orderable Group with Insoluble Word Problem.} \author{ V.\,V.~Bludov and A.\,M.\,W.~Glass } \small{\date{\today}} \maketitle { \def\thefootnote{} \footnote{2010 AMS Classification: 20F10, 06F15, 20F60.


This research was supported by a grant from the Royal Society for which we are most grateful. It was done while the first author was a Distinguished Academic Visitor at Queens' College, Cambridge; he wishes to thank the College and DPMMS for their warm hospitality.} } {\large Dedicated to James Wiegold, in memoriam.}

\bigskip

\setcounter{theorem}{0}
\setcounter{footnote}{0}

\begin{abstract}

We construct a finitely presented (two-sided) totally orderable group with insoluble word problem.

\end{abstract}

\section{Introduction}

Let $G$ be a group with a total order $<$. Then $G$ is said to be a
{\it right-ordered group} (with respect to $<$) if $f<g$ implies
$fx<gx$ for all $x\in G$. We can similarly define a left-ordered
group. A group with a total order that is respected by both right
and left multiplication is called a {\it two-sided ordered group},
or an {\it ordered group} or o-{\it group}, for short. A group that
can made into a right-ordered group with respect to some total order
is called a {\it right-orderable} group, and a group that can be made into
an o-group with respect to some total order is called a ({\it
two-sided}) {\it orderable} group. In \cite{BGGS} and \cite{BG09},
we constructed finitely presented right-orderable groups with
insoluble word problem. These examples are not two-sided orderable.
Here we prove
\bigskip

\begin{theorem} \label{IWPO} There is a finitely presented (two-sided) orderable group with insoluble word problem. \end{theorem}

Our proof involves taking a finitely presented {\it right}-orderable
group $F/N$ with insoluble word problem and using it to obtain a
finitely presented (two-sided) orderable group with insoluble word
problem. We do this as follows.  Let $F$ be the free group on $m$
generators and $u_1,\dots,u_n$ be the generators of $N$ (as a {\it
normal} subgroup of $F$). Let $G_0$ be a semidirect product of $F$
by a free group on $2m$ generators that normalises $N$. Take an
HNN-extension $G_1$ of $G_0$ with stable letter $t$ that fixes each
element of $N$. Let $T_0$ be the normal closure of $\langle
t\rangle$ by $G_0$ (equivalently, in $G_1$). We use the right order on $F/N$ to give a
right order $\prec$ on $G_0$ and thence an order on the generators
of the free group $T_0$ by $t^{f_1} < t^{f_2}$ if and only if
$f_1\prec f_2$. We will give an ordering of basic commutators in a
free group and derive a $G_1$-invariant order on $T_0$. 
We embed $T_0$ into its topological completion $T_0^{\ast}$, choose an appropriate subgroup $T<T_0^{\ast}$, and form a semidirect product of $T$ by a direct product of $G_0$ and a free group $F(\bar y)$ on $2m$ generators.
This is our (two-sided) orderable group $G$.
We derive the infinite set of relations $t^{-1}ut=u$ for all $u\in
N$ from a {\it finite} set of defining relations of $G$ which
include $t^{-1}u_it=u_i$ ($i=1,\dots,n$). In our construction, $t$ commutes (in
$G$) with $x\in F$ if and only if $x\in N$. It follows that $G$ has insoluble
word problem.

\section{Preliminaries (Groups)}

We will write $\ell(w)$ for the {\it length} of a reduced word $w$ in a free group. \medskip

For general information about HNN-extensions, see \cite{LS}, Chapter IV, Section 2. We summarise what we require. \medskip

Let $G$ be a group with isomorphic subgroups $A$ and $B$, say $\varphi: A\cong B$. The {\it HNN-extension} of $G$ relative to $A$, $B$, $\varphi$ is the group $$ G^{\#}=\langle G,\,t;\,t^{-1}at=a^\varphi\; (a\in A)\rangle.$$

{\bf Definition} (\cite{LS}, page 181) A sequence $g_0,\,t^{\varepsilon_1},\,g_1,\,\dots,\, t^{\varepsilon_n},\,g_n$ ($n\geq 0$, $\varepsilon_i\in \{ \pm 1\},~i=1,\dots,n$) is said to be {\it reduced} if there is no consecutive subsequence $t^{-1},\, g_i,\, t$ with $g_i\in A$ or $t,\, g_i,\, t^{-1}$ with $g_i\in B$. \bigskip

{\bf Britton's Lemma} (\cite{LS}, page 181) {\it If the sequence $g_0,\,t^{\varepsilon_1},\,g_1,\,\dots,\, t^{\varepsilon_n},\,g_n$ is reduced and $n\geq 1$, then $g_0t^{\varepsilon_1}g_1\cdots t^{\varepsilon_n}g_n\neq 1$ in $G^{\#}$. } \bigskip

If $C$ is any subgroup of a group $G$, we always let $1$ be the right or left coset representative of $C$. \bigskip

{\bf Definition.} (\cite{LS}, page 181) A {\it normal form} is a sequence

 $$g_0,\,t^{\varepsilon_1},\,g_1,\,\dots,\, t^{\varepsilon_n},\,g_n\;\;\; (n\geq 0),\;\;\;\;\;\;\;\hbox{where}$$

(i) $g_0$ is an arbitrary element of $G$,

(ii) if $\varepsilon_i=-1$, then $g_i$ is a representative of a right coset of $A$ in $G$,

(iii) if $\varepsilon_i=+1$, then $g_i$ is a representative of a right coset of $B$ in $G$, and

(iv) there is no consecutive subsequence $t^{\varepsilon}$, $1$, $t^{-\varepsilon}$ with $\varepsilon\in \{ \pm 1\}$. \bigskip

{\bf Theorem 2.1} (\cite{LS}, page 182) {\it Let $G^{\#} = \langle G,\,t;\,t^{-1}at=a^\varphi\; (a\in A)\rangle$ be an HNN-extension of $G$. Then \medskip

{\rm (I)} The group $G$ is embedded in $G^{\#}$ by the map $g\mapsto g$. If $g_0t^{\varepsilon_1}g_1\cdots t^{\varepsilon_n}g_n = 1$ in $G^{\#}$ where $n\geq 1$, then $g_0,\,t^{\varepsilon_1},\,g_1,\,\dots,\, t^{\varepsilon_n},\,g_n$ is not reduced. \medskip

{\rm (II)} Every $w\in G^{\#}$ has a unique representation as $w=g_0t^{\varepsilon_1}g_1\cdots t^{\varepsilon_n}g_n$ where $g_0,\,t^{\varepsilon_1},\,g_1,\,\dots,\, t^{\varepsilon_n},\,g_n$ is a normal form.} \medskip

The only background needed about basic commutators can be found in \cite{Ha}, Section 11.1.
\medskip

We will use the following terminology throughout. If $X$ and $Y$ are totally ordered sets, then we can totally order the set $X \times Y$ by $$(x_1,y_1)<(x_2,y_2)\quad \hbox{iff}\quad x_1<x_2\;\hbox{ or }\; (x_1=x_2\;\&\;y_1<y_2).$$ This is called the {\it lexicographic product} of $X$ and $Y$. If $X$ and $Y$ are right-ordered groups, then so is $X \times Y$ under this ordering; if $X$ and $Y$ are o-groups, then so is $X \times Y$ under this ordering.

\section{A first part of construction: the group $G_1$.}

Let $\hat H$ be any finitely presented right-orderable group with insoluble word problem; say, $\hat H:=\langle x_1,\dots,x_m;\, u_1=\dots=u_n=1\rangle$. (For the existence of such $\hat H$, see \cite{BGGS} or \cite{BG09}.) Let $F(\bar x):=F(x_1,\dots,x_m)$ denote the free group on free generators $x_1,\dots,x_m$ and write $\hat H=F(\bar x)/N$, where $N$ is the normal subgroup of $F(\bar x)$ generated by $u_1,\dots,u_n$.  \medskip

Let $G_0$ be generated by $x_1,\dots,x_m,\; b_1,\dots,b_{2m}$ and have defining relations

\begin{equation} \label{d5} 
x_i^{b_j}=x_i^{x_j},\quad \;\;\;\; x_i^{b_{m+j}}=x_i^{x_j^{-1}}\qquad (i,j=1,\dots,m). 
\end{equation}

\noindent So $G_0$ is a semidirect product of the free group $F(\bar x)$ by the free group $F(\bar b)$ with free generators $b_1,\dots,b_{2m}$, and $N$ is normalised by $F(\bar b)$ in $G_0$. \medskip

Let $G_1$ be generated by $x_1,\dots,x_m,\; b_1,\dots,b_{2m},$ and the extra generator $t$ and have defining relations (\ref{d5}) and

\begin{equation}
\label{d6}
[t,u_j^g]=1\;\;\;\;\;\;(j=1,\dots,n;\;~g\in F(\bar x)).
\end{equation}
\medskip

\noindent So $G_1$ is an HNN-extension of $G_0$. By Theorem 2.1, $G_0$ (and hence $F(\bar x)$ and $F(\bar b)$) can be embedded in $G_1$ in the natural way.
We will regard $N$, $F(\bar x)$ and $F(\bar b)$ as subgroups of $G_1$.
By Britton's Lemma,

\begin{lemma}
\label{L2}
If $w\in F(\bar x)$, then
$\;[t,w]=1\;\hbox{in}\;G_1\; \hbox{if and only if} \;\; w\in N.$
\end{lemma}

Let $T_0$ be the normal closure of $\langle t\rangle$ in $G_1$ (so $T_0=\langle t\rangle^{G_0}$).
\medskip

For each $\hat h\in \hat H=F(\bar x)/N$, choose $h\in F(\bar x)$ a preimage of $\hat h$. We will choose $1$ to be the preimage of $\hat 1$. So $Nh=hN=\hat h$ and $H:=\{ h\in F(\bar x)\mid \hat h\in \hat H\}$ is a transversal for $N$ in $F(\bar x)$. \medskip

For $f_1,f_2\in F(\bar x)$ and $v_1(\bar b),v_2(\bar b)\in F(\bar b)$, we have $$t^{v_1(\bar b)f_1}=t^{v_2(\bar b)f_2}\quad \hbox{iff}\quad (v_1(\bar b)=v_2(\bar b)\;\,\&\;\,\hat f_1=\hat f_2).$$ Hence, using the normal form for elements of an HNN-extension or Britton's Lemma (see above), one immediately obtains

\begin{lemma}\label{L1} 
$T_0$ is a free group with free generators $t^{vh}$
$(v\in F(\bar b);\;h\in H)$.
\end{lemma}

\section{Preliminaries (Orderability)}

 Throughout the rest of the paper, we will use $<$ for a two-sided total order and $\prec$ for a right total order on a group. \medskip

If $G$ is a right-ordered group, let $G_+:=\{ g\in G: g\succ 1\}$, the set of {\it strictly positive elements of} $G$. Note that $G_+$ is a subsemigroup of $G$ and $\{ G_+, G_+^{-1},\{ 1\}\}$ is a partition of $G$. If $P$ is a subsemigroup of $G$ and $\{ P, P^{-1},\{ 1\}\}$ is a partition of $G$, then $G$ can be right ordered by setting $f\prec g$ if and only if $gf^{-1}\in P$. If $G$ is an o-group, then $G_+$ is a normal subsemigroup of $G$ with the above properties, and any normal subsemigroup $P$ with these properties can be used to make $G$ an o-group.\medskip

Let $G$ be an ordered group and $C$ be a subgroup of $G$.
We say that $C$ is {\it convex in} $G$ if $g\in C$ whenever there are $c_1,c_2\in C$ with  $c_1<g<c_2$. The set of convex subgroups of $G$ form a totally ordered set under inclusion (\cite{G99}, Lemma 3.1.2). For a subgroup $H$ of an ordered group $G$, we denote by $con(H)$ the {\it convexification} of $H$ in $G$; $con(H)$ is the smallest convex subgroup of $G$ which contains $H$, so $con(H)$ is equal to the intersection of all convex subgroups of $G$ containing $H$. In two-sided ordered groups, $con(H)$ contains those and only those elements $g\in G$ for which there  are $h\in H$ such that $h^{-1}\leq g\leq h$.  
\medskip  

We call an o-group {\it archimedean} if (whenever $g\in G$ and $f^k<g$ for all $k\in \Z$,
then $f=1$). Every archimedean o-group is abelian and is
isomorphic to a subgroup of the additive group of real numbers
equipped with the usual ordering (H\"{o}lder's Theorem; see, {\it
e.g.}, \cite{G99}, Theorem 4.A).
\medskip

If $G$ is an o-group and $g,h\in G_+$, we will say that $g$ and $h$ are {\it archimedean equivalent} if the convex subgroups that they generate are equal; that is, there are $r,s\in \Z_+$ such that $g < h^r$ and $h<g^s$. Thus an o-group $G$ is archimedean if and only if all elements of $G_+$ are archimedean equivalent.
\medskip

If $G$ is an o-group and $g,h\in G_+$ with $g^k<h$ for all $k\in \Z_+$, we will write $g\ll h$.
\medskip

If $H$ is an o-group and the group $G$ acts on $H$, then we say that
the order on $H$ is $G$-{\it invariant} if $h>1$ implies $h^g>1$ for
all $h\in H$ and $g\in G$, and a subgroup $K$ of $H$ is $G$-{\it
invariant} if $K^g=K$ for all $g\in G$.
\medskip

We will also need that for any o-group $G$, the topological completion $G^\ast$ of $G$ under the order
topology is an o-group (see \cite{Ban}). The o-group $G^\ast$ can also be realised
as the sequential completion of $G$; that is, every sequence that is
a left and right Cauchy sequence is convergent, and all elements of $G^\ast$ are limits of Cauchy sequences in $G$ (see \cite{CH63}).
\medskip

Let $F$ be a free group. The standard way to order $F$ is to use the
lower central series $\gamma_k(F)$ ($k\in \Z_+$). The key is that
$\gamma_k(F)/\gamma_{k+1}(F)$ is a free abelian group and so can be
made into an o-group ($k\in \Z_+$). Since $\bigcap_{k\in \Z_+}
\gamma_k(F) =\{ 1\}$, we can produce a two-sided order on $F$ as
follows: Make each $\gamma_k(F)/\gamma_{k+1}(F)$ an o-group. Let
$f,g\in F$ with $f\neq g$. Let $s\in \Z_+$ be such that $gf^{-1}\in
\gamma_{s}(F) \setminus\gamma_{s+1}(F)$. Then $F$ is an o-group if
we define $f<g$ if and only if $\gamma_{s+1}(F)1<\gamma_{s+1}(F)gf^{-1}$
in $\gamma_{s}(F)/\gamma_{s+1}(F).$ We call any order constructed in
this way a {\it standard central order}; so standard central orders
on free groups depend only on the orders defined on the set of
abelian groups $\{\gamma_k(F)/\gamma_{k+1}(F)\mid k\in \Z_+\}$. For
this and further background, also see \cite{F}, 
\cite{G99} or \cite{KK72}.
\medskip

\section{Ordering $G_1$}\label{OrdG1}

As noted in the previous section, the free groups $F(\bar x)$ and $F(\bar b)$ are two-sided orderable groups. By (\ref{d5}), $F(\bar b)$ acts by conjugation on $F(\bar x)$. Hence $G_0$ is an o-group if we define 
\begin{equation}\label{ord0}
u(\bar b)v(\bar x)>1\quad \hbox{iff}\quad u(\bar b)>1 \;\;\;\hbox{or}\;\;\; (u(\bar b)=1\;\,\&\;\,v(\bar x)>1).
\end{equation}
Thus $G_1/T_0\cong G_0$ is a two-sided ordered group. Now $G_1$ is an ordered group with convex subgroup $T_0$ \underline{{\it if}} we can define a $G_1$-invariant (two-sided) order on $T_0$ (define $g>1$ in $G_1$ if and only if either $gT_0>T_0$ in $G_1/T_0$ or $g\geq 1$ in $T_0$). So it remains to construct a $G_1$-invariant order on $T_0$. Note that the construction will not be effective (it doesn't need to be). Really we use only one non-effective step which is confined to the existence of a right order on $\hat H$.
\medskip

To construct such an order on $T_0$, we first put a {\it right}
order $\prec$ on $G_0$ as follows. Let $\prec$ be the right total
order on $\hat H$. Since $N$ is a subgroup of a free group and hence
is free, it can be made into an o-group as described in the previous
section using any standard central ordering. Define a right order on $F(\bar x)$ by:
$$f\succ 1\quad \hbox{iff}\quad Nf\succ N\;\;\;\hbox{or}\;\;\;f\in
N_+.$$ Next, as described in the previous section, we can  put a
standard central order on $F(\bar b)$ so that $b_i <1$
($i=1,\dots,2m$) and the order on the free abelian group $F(\bar
b)/\ga_2(F(\bar b))$ is archimedean. Define 
$$v(\bar b)h\succ 1\quad \hbox{iff}\quad v(\bar b)\in F(\bar b)_+\;\;\hbox{or}\;\; (v(\bar
b)=1\;\,\&\,\;h\in F(\bar x)_+).$$ 
This is a well-defined {\it right} order on the group $G_0$. Let
$$\Lambda:=\{ v(\bar b)h\in G_0\mid v(\bar b)\in F(\bar b),\;h\in
H\}$$ with the inherited order.
\medskip

By Lemma \ref{L1}, $\{ t^g\mid g\in \Lambda\}$ is a set of free generators for $T_0$; it inherits a total order $\ll$ given by: 
$$t^f\ll t^g\quad \hbox{iff}\quad f\prec g\quad (f,g\in \Lambda).$$ 
So for each $a\in G_0$, $t^{fa}\ll t^{ga}\;$ if and only if $\;t^f\ll t^g$. \medskip

For each $g\in \Lambda$, consider the normal subgroups of $T_0$ $$C_g:=\langle t^f\mid f\in \Lambda,\;\; f\prec g\rangle^{T_0}\quad \hbox{and}\quad C(g):=\langle t^f\mid f\in \Lambda,\;\; f\preceq g\rangle^{T_0}.$$
We will use this set of normal subgroups of $T_0$ to define a two-sided order on $T_0$ that is $G_1$-invariant (see Lemma \ref{0}). 
\medskip

For each $g\in \Lambda$, let 
$$K_g:=\langle t^f\mid f\in \Lambda,\;\; f\succeq g\rangle\quad \hbox{and}\quad K(g):=\langle t^f\mid f\in \Lambda,\;\; f\succ g\rangle.$$
By Lemma \ref{L1}, $K_g$ and $K(g)$ are free groups on the indicated generators ($g\in \Lambda$).

\begin{lemma}
\label{free0}
Let $g\in \Lambda$.
\medskip

{\rm (i)} $[t^f,t^{f'}]\in C(g)$ for each $f,f'\in G_0$ with $f\preceq g$.
\medskip

{\rm (ii)} $T_0/C(g)\cong K(g)$, $T_0/C_g\cong K_g$; and $T_0$ is isomorphic to the semidirect product $C(g)\rtimes K(g)$ as well as to the semidirect product $C_g\rtimes K_g$. \medskip

{\rm (iii)} $\;C(g)^f=C(gf)$ and $C_g^f=C_{gf}$ for each $f\in G_0$.
\medskip

\noindent Moreover,
\medskip

{\rm (iv)} $\;T_0 = \bigcup_{g\in \Lambda} C(g) = \bigcup_{g\in \Lambda} C_g,$ and
\medskip

{\rm (v)} $\;\bigcap_{g\in \Lambda} C(g) = \bigcap_{g\in \Lambda} C_g =\{ 1\}.$
\end{lemma}

{\it Proof}: This is immediate by the definitions and our choice of ordering of the generators of $T_0$. \pfend

\begin{lemma} \label{free11} For any $w\in T_0\setminus \{ 1\}$, there is a unique $g\in\Lambda$ such that $w\in C(g)\setminus C_g$ \end{lemma}

{\it Proof}: Any $w\in T_0\setminus \{ 1\}$ can be written as a reduced word in some $t^{f_1},\dots, t^{f_n}$ with $f_1\prec \dots \prec f_n$ in $\Lambda$. Clearly, $w\in C(f_n)$. Choose the minimal $k\in\{1,\dots,n\}$ such that $w\in C(f_k)$. Then $w\not\in C_{f_k}$ (otherwise, $w\in \langle t^f\mid f\in \Lambda,\;\; f\prec f_k\rangle^{T_0}$; so $w\in C(f_{k-1})$, a contradiction). \pfend

We now show how to construct a $G_1$-invariant order on $T_0$. We first define a lexicographic order on the abelian groups $\gamma_k(T_0)/\gamma_{k+1}(T_0)$ ($ k\in\Z_+ $). For this we associate with any basic commutator $w(t^{f_1},t^{f_2},\dots,t^{f_s})$, the monomial $s(w)=f_1^{n_1}f_2^{n_2}\dots f_s^{n_s}$, where $n_i$ is the number of occurrences of $t^{f_i}$ in $w$ and $f_1\succ f_2\succ\dots \succ f_s$.
For example, if $w=[t^{f_1},[t^{f_1},t^{f_5}], [t^{f_2},t^{f_3},t^{f_5},t^{f_5}]]$ and $f_1\succ \dots \succ f_5$, then $s(w)=f_1^2f_2f_3f_5^3$. We put the lexicographical order on the set of such monomials. For basic commutators $w_1,w_2\in \gamma_k(T_0)\setminus\gamma_{k+1}(T_0)$, define $$ w_1\ll w_2\quad\mbox{if} \quad s(w_1)<s(w_2);$$ and fix $w_1\ll w_2$ arbitrarily if $s(w_1)=s(w_2)$. Now let $\leq_0$ be the standard central order on $T_0$ built from the total orders on $\gamma_k(T_0)/\gamma_{k+1}(T_0)$ ($k\in \Z_+$) described above. The order $\leq_0$ induces a total order on the subgroup $K_1$. The isomorphism given in Lemma~\ref{free0}(ii) induces a total order on the group $T_0/C_1$. This, in turn, induces an order on its subgroup $C(1)/C_1$. We also denote this order by $\leq_0$. Let $w$ be a non-trivial element from $T_0$. By Lemma~\ref{free11}, there is a unique $g\in\Lambda$ such that $w\in C(g)\setminus C_g$. Define 
\begin{equation}\label{ord1}
  1<_1 w \quad\mbox{iff}\quad C_1<_0 C_1w^{g^{-1}}\mbox{ in } C(1)/C_1.
\end{equation}

\begin{lemma}\label{Cg}
Subgroups $C_g$ and $C(g)$ are convex in $G_1$ under order $\leq_1$ for all $g\in \Lambda$.
\end{lemma}

{\it Proof}: Since $T_0$ is convex in $G_1$, it is sufficient to prove that $C_g$ and $C(g)$ are convex in $T_0$. Let $1 <_1 u <_1 h$, $h\in C(g)$. This gives $1 <_1 hu^{-1}$. By Lemma~\ref{free11}, we find $f\in\Lambda$ such that $u\in G(f)\setminus G_f$. If $u\not\in C(g)$, then $g\prec f$ and so $C(g)\leq C_f$.  
Thus $u^{-1}\in G(f)\setminus G_f$ and $hu^{-1}\in G(f)\setminus G_f$. By (\ref{ord1}), $C_1<_0 C_1(hu^{-1})^{f^{-1}} = C_1h^{f^{-1}}u^{-f^{-1}} = C_1u^{-f^{-1}}$ and so $1 <_1 u^{-1}$, a contradiction. Therefore $u\in C(g)$ and $C(g)$ is convex. The convexity of $C_g$ follows and is left to the reader.\pfend

\begin{lemma}\label{0} 
The order $\leq_1$ defined in $(\ref{ord1})$ is a
$G_1$-invariant two-sided order on $T_0$.
\end{lemma}

{\it Proof}:
 Assume that $1<_1 w_1,w_2$.
By Lemma \ref{free11}, there are $g_1,g_2\in \Lambda$ such that
$w_1\in C(g_1)\setminus C_{g_1}$ and $w_2\in C(g_2)\setminus
C_{g_2}$. If $g_1\prec g_2$, then $w_1w_2\in C(g_2)\setminus
C_{g_2}$ and $C_1 <_0 C_1w_2^{g_2^{-1}}= C_1(w_1w_2)^{g_2^{-1}}$.
Hence $1<_1 w_1w_2$. Similarly, $1 <_1 w_1w_2$ if $g_2\prec g_1$. If
$g_1=g_2$, then $C_1 <_0 C_1w_1^{g_2^{-1}}$ and $C_1 <_0
C_1w_2^{g_2^{-1}}$; so $C_1 <_0 C_1(w_1w_2)^{g_2^{-1}}$, whence
$1<_1 w_1w_2$.
\medskip

Let $w\in T_0$ and $1<_1 w$. By Lemma \ref{free11}, there is $g\in\Lambda$ such that $w\in C(g)\setminus C_g$ and $C_1<_0 C_1w^{g^{-1}}$. Since $C_g$ is a normal subgroup of $T_0$, we get $w^v \in C(g)\setminus C_g$ for any $v\in T_0$. To prove that $<_1$ is a two-sided total order on $T_0$, we must show that $C_1 <_0 C_1(w^v)^{g^{-1}}$. This immediately follows from the definition of a standard central order on $T_0$ because $\gamma_{k+1}(T_0)w=\gamma_{k+1}(T_0)w^v$ provided that $w\in\gamma_k(T_0)\setminus\gamma_{k+1}(T_0)$. Thus the order $\leq_1$ is $T_0$-invariant; {\it i.e.}, $\leq_1$ is a two-sided order on $T_0$. \medskip

Let $f\in\Lambda$. Now $w^f\in C(gf)\setminus C_{gf}$ and $(w^f)^{(gf)^{-1}}=w^{g^{-1}}$. Hence $1<_1w^f$. Thus the order $\leq_1$ is $G_1$-invariant.  \pfend

\section{The second part of construction: the group $G$.}

Denote by $T_0^{\ast}$ the topological (sequential) completion of $T_0$ under the interval topology induced by the order on $T_0$ defined in Section \ref{OrdG1}. We regard $T_0^{\ast}$ as an o-group with the order $\leq_{\ast}$ extending the initial order $\leq_1$ of the group $T_0$. Since convex subgroups $C(g)$ are normal in $T_0$ (by Lemma~\ref{Cg}) and their intersection is trivial (by Lemma~\ref{free0}~(v)), the right and the left topological spaces coincide on $T_0$. In this case, any right Cauchy sequence is a left Cauchy sequence and {\it vice versa} and we can simply write ``Cauchy sequence" without ambiguity.
Conjugation by $g\in G_1$ preserves the order on $T_0$, so maps open intervals on open intervals; hence
conjugation is a continuous operation. This allows us to define an order-preserving action by $G_0<G_1$ on $T_0^{\ast}$ coinciding with conjugation on $T_0$.
Let $\{ c_k\mid k\in \Z_+\}$ be a Cauchy sequence in $T_0$. Then
\begin{equation}\label{ck}
(\lim_{k\rightarrow \infty} c_k)^g=\lim_{k\rightarrow \infty} c_k^g \quad (\lim_{k\rightarrow \infty} c_k \in T_0^{\ast},~c_k\in T_0,~g\in G_0).
\end{equation}
We define now $2m$ order-preserving automorphism $y_1,\dots,y_{2m}$ on $T_0^{\ast}$.
First, let
\begin{equation}\label{top}
(t^g)^{y_i}=[b_i,t]^g\quad (g\in \Lambda;\;i=1,\dots,2m).
\end{equation}
By Lemma \ref{L1}, the set $\{ t^g\mid g\in \Lambda\}$ freely generates $T_0$ so each $y_i$ uniquely defines an endomorphism $y_i: T_0\rightarrow T_0$, 
$$w(t^{g_1},\dots,t^{g_k})^{y_i}=w(t^{g_1y_i},\dots,t^{g_ky_i})\quad (w(t^{g_1},\dots,t^{g_k})\in T_0,~i=1,\dots,2m).$$ 
It follows from (\ref{top}) that 
\begin{equation}\label{comy}
w^{y_ig}=w^{gy_i}\quad \mbox{for all } w\in T_0,~g\in \Lambda;\;i=1,\dots,2m.
\end{equation}
\begin{lemma} \label{ordy} The endomorphisms $y_1,\dots,y_{2m}$ of $T_0$ preserve the order $\leq_1$ on $T_0$; moreover $w^{y_i}$ is archimedean equivalent to $w$ for all $w\in T_0$, $i=1,\dots, 2m$ . \end{lemma}
{\it Proof}: 
Consider a basic commutator $d\in\gamma_k(T_0)/\gamma_{k+1}(T_0)$, $k\in\Z_+$ and assume $1<_0 d$. Applying (\ref{top}) to $d=[t^{f_1},\dots,[\dots,\dots],t^{f_s}]$, we obtain (modulo $\ga_{k+1}(T_0)$) that
\begin{equation}\label{dy}
d^{y_i}=[t^{-b_if_1}t^{f_1},\dots,[\dots,\dots],t^{-b_if_s}t^{f_s}]=d\cdot d_1\dots d_{p},
\end{equation} 
where $s(d_1),\dots,s(d_{p})<s(d)$. By the definition of order $\leq_0$, all $d_1,\dots, d_p\ll d$ and so $d^{y_i}>_0 1$ and $d^{y_i}$ is archimedean equivalent to $d$ ($i=1,\dots,2m$). Now we repeat arguments from the proof of Lemma~\ref{0}.
Let $k\in\Z_+$ be such that $w^{g^{-1}}\in \gamma_k(T_0)\setminus\gamma_{k+1}(T_0)$. So there are basic commutators $c_1\gg \dots \gg c_\ell$ in $\gamma_k(T_0)$ and $r_1,\dots,r_\ell\in \Z\setminus \{ 0\}$ such that $$w^{g^{-1}}=c_1^{r_1}\dots c_\ell^{r_\ell}w',$$ with $w'\in \ga_{k+1}(T_0)$. Since $1<_1 w$, we have $r_1\in \Z_+$ by the definition of the order $\leq_1$. Congugating $w^{g^{-1}}$ by $y_i$ and using  (\ref{comy}) gives $c_1^{r_1y_i}>_1 1$ and 
$c_1^{y_i}$ is archimedean equivalent to $c_1$ by (\ref{dy}). Thus the order $\leq_1$ is invariant under each of the endomorphisms $y_1,\dots,y_{2m}$ and $w^{y_i}$ is archimedean equivalent to $w$ for all $w\in T_0$. \pfend
\begin{lemma} \label{L9} The endomorphisms $y_1,\dots,y_{2m}$ of $T_0$ extend to order-preserving automorphisms of $T_0^*$. \end{lemma}
{\it Proof}:
First we show that endomorphism $y_i$ maps a Cauchy sequence on a Cauchy sequence. Indeed, if $\{c_k\}_{k\in\Z_+}$ is a Cauchy sequence, then for any $g\in\Lambda$ there is $k_g\in\Z_+$ such that $c_{k_g}c^{-1}_{k_g+s}\in C(g)$ for all $s\in\Z_+$. By Lemma~\ref{ordy}, $c_{k}$ is archimedean equivalent to $c^{y_i}_{k}$ and so $c_{k_g}^{y_i}c^{-y_i}_{k_g+s}\in C(g)$. Thus $\{c_k^{y_i}\}_{k\in\Z_+}$ is a Cauchy sequence.
For each $g\in G_0$ and $i=1,\dots,2m$ consider sequences $\{c_k(g,i)\}_{k\in\Z_+}$ defined by
\begin{equation}\label{ckg}
c_0(g,i):=t^g,\quad c_k(g,i):=t^{b_i^kg}c_{k-1}(g,i),\qquad k\in \Z_+.
\end{equation}
The ordering on $G_1$ ensures that $\{c_k(g,i)\}_{k\in\Z_+}$ is a Cauchy sequence in $T_0$ and so
$\tilde c(g,i):= \lim_{k\rightarrow \infty} c_k(g,i) \in T_0^{\ast}$ for every $g\in G_0$ and $i=1,\dots,2m$. By routine verification, $$(c_k(g,i))^{y_i}=t^{-b_i^{k+1}g}t^g\quad(k\in \Z_+).$$ 
By the sequential completeness of $T_0^{\ast}$ we have 
$$(\tilde c(g,i))^{y_i}=\lim_{k\rightarrow \infty} (c_k(g,i))^{y_i}=t^g.$$ 
It follows now that each element $w=w(t^{g_1},\dots,t^{g_k})\in T_0$ has a unique preimage $w^{y_i^{-1}}=w(t^{g_1y_i^{-1}},\dots, t^{g_ky_i^{-1}})\in T_0^{\ast}$ ($i=1,\dots,2m$). Given a Cauchy sequence $\{w_k\}_{k\in\Z_+}$, the sequence $\{w_k^{y_i^{-1}}\}_{k\in\Z_+}$ is also a Cauchy sequence. Thus $(\lim_{k\rightarrow \infty} w_k)^{y_i^{-1}}=\lim_{k\rightarrow \infty} w_k^{y_i^{-1}}\in T_0^{\ast}$. Hence each $\tilde w\in T_0^{\ast}$ has a unique preimage $\tilde w^{y_i^{-1}}$. So $y_i$ extends to an order-preserving automorphism of $T_0^*$ ($i=1,\dots,2m$). \pfend

Let $F(\bar y)$ be a free group on free generators $y_1,\dots,y_{2m}$. By Lemma~\ref{L9}, $F(\bar y)$ acts on $T_0^{\ast}$ by order-preserving automorphisms whose action on $T_0$ is given by (\ref{top}). We also defined the order-preserving action of $G_0$ on $T_0^{\ast}$ above in (\ref{ck}). Since these two actions are continuous and commute on $T_0$, they commute on $T_0^{\ast}$, and we can form the semidirect product $T_0^{\ast}\rtimes(G_0\times F(\bar y))$. The group $T_0^{\ast}\rtimes(G_0\times F(\bar y))$ is orderable because $T_0^{\ast}$ and $G_0\times F(\bar y)$ are orderable and $G_0\times F(\bar y)$ acts an $T_0^{\ast}$ by order-preserving automorphisms (with respect to $\leq_{\ast}$). 
Define
\begin{equation}\label{G}
G:=\langle t, x_1,\dots,x_m, b_1,\dots,b_{2m}, y_1,\dots,y_{2m}\rangle<T_0^{\ast}\rtimes(G_0\times F(\bar y)),
\end{equation}
and
\begin{equation}\label{T}
T:=\langle t \rangle^G<T_0^{\ast}.
\end{equation}
Thus we obtain

\begin{proposition}
\label{prop1}
The group $G$ is orderable.
\end{proposition}

It follows from the defintion (\ref{G}) the the group $G$ is generated by $5m+1$ elements: $t, x_1,\dots,x_m, b_1,\dots,b_{2m}, y_1,\dots,y_{2m}$ and $G$ contains the subgroup $\langle t, x_1,\dots,x_m, b_1,\dots,b_{2m}\rangle$ isomorphic with $G_1$. Hence $G$ satisfies relations (\ref{d5}) and (\ref{d6}) which hold in $G_1$. By the construction, $G$ satisfies relations (\ref{top}) and (\ref{comy}). In addition, $G$ satisfies the relations
\begin{equation}
\label{d1}
[x_i,y_j]=1\;\;\;\;\;\;(i=1,\dots,m;\;\;j=1,\dots,2m),
\end{equation}
and
\begin{equation}
\label{d2}
[b_i,y_j]=1\;\;\;\;\;\;(i,j=1,\dots,2m)
\end{equation}
which follow from (\ref{G}). Finally we have that the group $G$ satisfies relations (\ref{d5}), (\ref{d6}), (\ref{top}), (\ref{comy}), (\ref{d1}), and (\ref{d2}). We extract a finite subset from this set of relations, namely the relation (\ref{d5}), (\ref{d1}), (\ref{d2}) and relations
\begin{equation}
\label{d3} t^{y_i}=[b_i,t] \;\;\;\;\;\;(i=1,\dots,2m),
\end{equation}
and
\begin{equation}
\label{d4}
[t,u_j]=1\;\;\;\;\;\;(j=1,\dots,n).
\end{equation}

\medskip

\begin{lemma}\label{RelG} 
The two sets of relations relations $(\ref{d5})$, $(\ref{d6})$, $(\ref{top})$, $(\ref{comy})$, $(\ref{d1})$, $(\ref{d2})$  and $(\ref{d5})$, $(\ref{d1})$ -- $(\ref{d4})$ are equivalent.
\end{lemma}

{\it Proof:} First we show that relations (\ref{d5}), (\ref{d6}), (\ref{top}), (\ref{comy}), (\ref{d1}), (\ref{d2}) follow from relations $(\ref{d5})$ and $(\ref{d1})$ --- $(\ref{d4})$.

Conjugating $[t,u_j]=1$ by $y_i$ and then by $b_i^{-1}$
(and by $y_{i+m}$ and then by $b_{i+m}^{-1}$) and using (\ref{d5}),
(\ref{d1}) --- (\ref{d3}), we obtain $[t,u_j^{x_i^{\pm 1}}]=1$ ($j=1,\dots, n;\;i=1,\dots,m$). An easy induction now gives that $[t,u_j^{w(\bar x)}]=1$ for all $j\in \{ 1,\dots,n\}$ and
$w(\bar x)\in F(\bar x)$. Hence the relations (\ref{d6}) hold in $G$. Relations (\ref{top}) follow from (\ref{d3}) and (\ref{d1}), (\ref{d2}). Relations (\ref{comy}) follow from (\ref{d1}), (\ref{d2}). 

Now, relations (\ref{d3}), (\ref{d4}) are the partial cases of relations (\ref{top}) and (\ref{d6}) respectively. The lemma follows. \pfend 

Note that
$G\cong T\rtimes (G_0\times F(\bar y))$ by construction and so
\begin{equation}\label{GT}
G/T\cong G_0\times F(\bar y)
\end{equation}
However the isomorphsm (\ref{GT}) follows immediately from $(\ref{d5})$, $(\ref{d1})$ -- $(\ref{d4})$  

\medskip

Let $w\in F(\bar x)$. Since $[t,w]=1$ in $G$ if and only if
$w\in N$ (by Lemma~\ref{L2}), it follows that 

\begin{proposition}
\label{prop2}
The group $G$ has insoluble word problem.
\end{proposition}

\medskip
Therefore, the rest of the paper is devoted to showing that $G$ is finitely presented. By Lemma~\ref{RelG}, it is enough to prove that relations (\ref{d5}), (\ref{d6}), (\ref{top}), (\ref{comy}), (\ref{d1}), (\ref{d2}) completely define the group $G$. To achieve this, in the next section we construct a generating set for the subgroup $T$. We complete this section by considering convexifications (convex closures) of the subgroups $C_g$ and $C(g)$ in the group $G$. For $g\in\Lambda$, denote the convexifications of $C_g$ and $C(g)$ in $G$ by $D_g$ and $D(g)$, respectively. So

$$D_g=con(C_g)\ \mbox{ and }\ D(g)=con(C(g)).$$ 

\begin{lemma} \label{free1} For each $g\in \Lambda$, 
\medskip

{\rm (i)} for each $u\in T$ there exists $v\in T_0$ such that $D_g u= D_g v$ and $D(g) u= D(g) v$;
\medskip

{\rm (ii)} $D_g$ and $D(g)$ are $\langle T,F(\bar y)\rangle$-invariant;
\medskip;
\medskip

{\rm (iii)} $D(g)\bigcap G_1=C(g)$ and $D_g\bigcap G_1=C_g$;
\medskip

{\rm (iv)} $T/D_g\cong K_g$ and $T/D(g)\cong K(g)$;
\medskip

{\rm (v)} $T\cong D_g\rtimes K_g\cong D(g)\rtimes K(g)$. 
\medskip

{\rm (vi)} $D_g^f=D_{gf}$ and $D(g)^f=D(gf)$ for all $f\in G_0$.
\medskip

\noindent Moreover,
\medskip

{\rm (vii)} $\bigcup_{g\in \Lambda} D(g) = T.$
\medskip

\end{lemma}

{\it Proof}: Let $u\in T$, say $u= \lim_{k\rightarrow \infty} v_k$ for some Cauchy sequence $\{ v_k\}_{k\in \Z_+}$ in $T_0$. Then there exists $k=k(g)\in \Z_+$ such that $u\cdot v_{k(g)}^{-1}\in C_g$. Since $C_g\leq D_g$, we get $D_g u=D_g v_{k(g)}$ and (i) is proved for $D_g$. 
\medskip

Now $D_g=con(C_g)$ is normalised by $T_0$ since $C_g$ is normal in $T_0$. For $u\in T$, by (i) there is $v\in T_0$ such that $uv^{-1}\in D_g$. Thus $D_g^u=D_g^v=D_g$ and $D_g$ is normal in $T$. Since $w^y_i$ is archimedean equivalent to $w$ and $D_g$ is convex, it follows that $D_g$ is $F(\bar y)$-invariant. Hence (ii) holds for $D_g$.
\medskip

Let $u\in G_1$ with $h_1\leq u\leq h_2$ for some $h_1,h_2\in C_g$. Since $C_g$ is convex in $G_1$, we get $u\in C_g$. This gives (iii) for $D_g$.
\medskip

Applying (i), we have $T=D_gT=D_gT_0$. By (iii), $D_g\bigcap T_0=C_g$ and hence $T/D_g=D_gT_0/D_g\cong T_0/(D_g\bigcap T_0)=T_0/C_g$.
\medskip

By Lemma~\ref{free0}, $T_0/C_g\cong K_g$; so $T/D_g\cong K_g$. Since $D_g\bigcap K_g\leq T_0$, we have $D_g\bigcap T_0=C_g$. Moreover, $D_g\bigcap K_g=\{1\}$ since
$C_g\bigcap K_g=\{1\}$. This gives (v) for $D_g$.
\medskip

By Lemma~\ref{free0} (iii), we have $$D_g^f=con(C_g)^f=con(C_g^f)=con(C_{gf})=D_{gf}.$$

By Lemma~\ref{free0} (iv), $T_0=\bigcup_{g\in \Lambda} C_g \subset \bigcup_{g\in \Lambda} D_g$. Since $\bigcup_{g\in \Lambda} D_g$ is convex and each $w\in T$  is archimedean equivalent to some $u\in T_0$, we get $T\subseteq \bigcup_{g\in \Lambda} D_g$. This gives (vii). 
\medskip

Similarly, we obtain (i) -- (vii) for $D(g)$. 
Hence the lemma follows.\pfend

\section{A generating set for $T$.}

We will need two identities that follow immediately from (\ref{d2}) and (\ref{d3}):
for $i=1,\dots,2m$,

\begin{equation}
\label{y1}
t^{y_i^{-1}b_i}=t^{y_i^{-1}}t^{-1},
\end{equation}
and
\begin{equation}
\label{y2}
t^{y_i^{-1}b_i^{-1}}=t^{y_i^{-1}}t^{b_i^{-1}}.
\end{equation}

\medskip

\begin{lemma} \label{L5} $T$ has generators 
\begin{equation} \label{gen} t^{\alpha(\bar y)v(\bar b)h(\bar x)}, 
\end{equation} 
where $h(\bar x)\in H$, $v(\bar b)\in F(\bar b)$ and $\alpha(\bar y)$ is either empty or, for some $i\in \{ 1,\dots,2m\}$, $\alpha(\bar y)$ is a non-trivial element of $F(\bar y)$ that begins with $y_i^{-1}$ and $v(\bar b)$ does not begin with $b_i^{\pm 1}$. 
\end{lemma}

{\it Proof}: The elements of $T$ of the form
\begin{equation}\label{gen'}
t^{\alpha(\bar y)v(\bar b)h(\bar x)}\quad (\alpha(\bar y)\in F(\bar y),\, v(\bar b) \in F(\bar b),\, h(\bar x)\in F(\bar x))
\end{equation}
generate $T$. Since $t^{\alpha(\bar y)v(\bar b)f_1(\bar x)}=t^{\alpha(\bar y)v(\bar b)f_2(\bar x)}$ if $Nf_1=Nf_2$, we may assume that $h(\bar x)\in H$ in (\ref{gen'}). \medskip

We prove by induction on pairs of natural numbers $(\ell(\alpha),\ell(v))$ (ordered lexicographically) that $t^{\alpha(\bar y)v(\bar b)h}$ can be written as a product of conjugates of $t$ and $t^{-1}$ all of the form described in (\ref{gen}). \medskip

If $\ell(\alpha)=0$, it already has the desired form. If $\alpha(\bar y)$ begins with $y_i$ for some $i\in \{ 1,\dots, 2m\}$, $\alpha:=y_i\cdot \alpha'$ then $$t^{\alpha(\bar y)v(\bar b)h}=t^{y_i\alpha'(\bar y)v(\bar b)h}.$$ By (\ref{d3}) and (\ref{d2}), $$t^{y_i\alpha'(\bar y)v(\bar b)h}=
 t^{-b_i\alpha'(\bar y)v(\bar b)h}\cdot t^{\alpha'(\bar y)v(\bar b)h}=
 t^{-\alpha'(\bar y)b_iv(\bar b)h}\cdot t^{\alpha'(\bar y)v(\bar b)h}.$$
The two conjugators on the right-hand side have $\bar y$-length
$\ell(\alpha')<\ell(\alpha)$, so (by induction) each can be written
as a product of the desired form. Hence, so can $t^{\alpha(\bar
y)v(\bar b)h}$.
\medskip

\noindent We may therefore assume that $\alpha(\bar y)$ begins with
$y_i^{-1}$ for some $i\in \{ 1,\dots, 2m\}$; say, $\alpha(\bar
y):=y_i^{-1}\alpha'(\bar y)$. If $v(\bar b)$ does not begin with
$b_i$ or $b_i^{-1}$, then $t^{\alpha(\bar y)v(\bar b)h}$ has the
desired form.
\medskip

If, on the other hand,
$v=b_iv'$ with $\ell(v')<\ell(v)$, then by (\ref{d2}) and (\ref{y1})
$$t^{\alpha(\bar y)v(\bar b)h}=t^{y_i^{-1}b_i\alpha'(\bar y)v'(\bar b)h}=
t^{\alpha(\bar y)v'(\bar b)h}\cdot t^{-\alpha'(\bar y)v'(\bar b)h}$$
and $(\ell(\alpha'),\ell(v'))<(\ell(\alpha),\ell(v'))<
(\ell(\alpha),\ell(v))$. By induction, we get that $t^{\alpha(\bar
y)v(\bar b)h}$ is the product of of conjugates of $t$ and $t^{-1}$
all of the form given in (\ref{gen}). If $v$ begins with $b_i^{-1}$,
we can repeat the above argument using (\ref{y2}) instead of
(\ref{y1}). This completes the proof that the elements displayed in
(\ref{gen}) are indeed generators of $T$. \pfend

We will write {\it Gen} for the set of generators described in (\ref{gen}); in Corollary \ref{L6}, we will prove that they form a {\it free} generating set.
\medskip

We will frequently use that for each $f,g\in \Lambda$ and $i\in \{1,\dots,2m\}$, there is $k=k_{i,f,g}\in \Z_+$ such that $b_i^kf\prec g$. This is immediate from the definition of $\prec$.
\medskip

For each $g\in \Lambda$, let $A(g):=(T/D(g))/\ga_2(T/D(g))$, a free abelian group by Lemma~\ref{free1} (iv). 

\begin{remark} \label{det} {\rm Let $g\prec b_i^k$ and $a_{k,p}$ be the coefficient of $b_i^k$ in $t^{y_i^{-p}}D(g)$ in $A(g)$. Then $a_{k,1}=1$ and $a_{k,p+1}=\sum_{\ell\leq k} a_{\ell,p}$. If $k\in \N$, $r\geq 2$ and $p_1,\dots,p_r\in \Z_+$ are distinct with $g\prec b_i^{k+r}$, then this recursive formula immediately gives that the matrix $(a_{\ell,p})$ with $1\leq \ell,p\leq P$ has determinant $1$ where $P:=\max\{ p_s\mid s=1,\dots,r\}$.
Thus the $r$ columns $(a_{\ell,p_s})$ are linearly independent, whence, using the recursive formula again, the $r\times r$ matrix $(a_{k+j,p_s})$ has non-zero determinant.} \end{remark}

\begin{lemma} \label{key} Let $\alpha_1,\dots,\alpha_r\in F(\bar y)$ be distinct with each $\alpha_s$ either empty or beginning with some $y_{i_s}^{-1}$. Then there is $g\in \Lambda$ such that $$\{ t^{\alpha_1}D(g),\dots,t^{\alpha_r}D(g)\}\;\;\hbox{ is linearly independent in } \;A(g).$$ \end{lemma}

{\it Proof}: Induction on $r$.
\medskip

Let $r=1$. By Lemma \ref{free1}(I), if $g\prec 1$ and
$g\succ b_1,\dots,b_{2m}$, then $t^{\alpha_1}D(g)=tD(g)$.
The result follows at once in this case.
\medskip

Assume the result if $r<r_0$ and let $t^{\alpha_1},\dots,t^{\alpha_{r_0}}$ satisfy the hypotheses of the Lemma. \medskip

If some $\alpha_{s_0}$ is empty, let $P_0=\{ t^{\alpha_{s_0}}\}=\{ t\}$. For all other values of $s$, write $\alpha_s=y_{i_s}^{-1}\beta_s$ where $\beta_s$ does not begin with $y_{i_s}$. Let $P_i=\{ t^{\alpha_s}\mid i_s=i\}$ ($i=1,\dots,2m$). The non-empty sets among $P_0,\dots,P_{2m}$ partition $t^{\alpha_1},\dots,t^{\alpha_{r_0}}$. Now, in $A(g)$, we have that $t^{(\alpha_s)_{1,g}}$ is a product which has terms ending in $b_{i_s}$. By Lemma \ref{L1}, the set of elements in $P_i$ is linearly independent in $A(g)$ of the set of elements in all the remaining $P_{i'}$ ($i'\neq i$). So if there are distinct $i,j$ with $P_i$ and $P_{j}$ non-empty, the original set is linearly independent in $A(g)$ by the induction hypothesis. We may therefore assume that there is a unique $i\in \{ 1,\dots,2m\}$ such that $P_i\neq \emptyset$. \medskip

Let $\alpha_s(\bar y)$ begin with $y_i^{-n_{1,s}}$ ($n_{1,s}\in \Z_+$) for $s=1,\dots,r_0$; say $$\alpha_s(\bar y)=y_{i}^{-n_{1,s}}\dots\; y_{i_{q_s,s}}^{n_{q_s,s}}.$$ For each $n\in \Z_+$, let $Q_n:=\{ s\in \{ 1,\dots,r_0\}\mid n_{1,s}=n\}.$ By Remark \ref{det}, the set of elements $t^{\alpha_s}$ with $s\in Q_n$ is linearly independent of $\{ t^{\alpha_{s'}}\mid s'\not\in Q_{n}\}$. So the lemma follows by induction unless $n_{1,s}=n_{1,1}$ for all $s=1,\dots,r_0$. We therefore assume that $$\alpha_s(\bar y)=y_{i}^{-n_{1,1}}\dots\; y_{i_{q_s,s}}^{n_{q_s,s}}=:y_{i}^{-n_{1,1}}\delta_s,$$ where $\delta_s$ starts with $y_{j_s}^{\pm 1}$ with $j_s\neq i$ ($s=1,\dots,r_0$).\medskip

We prove the result for $r_0$ by induction on $\ell=\max\{
\ell(\alpha_1),\dots,\ell(\alpha_{r_0})\}$. Now if $g\prec
b_i^{n_{1,1}}$, in $A(g)$ we have a linear combination of elements
all of the form $t^{\delta_1b_i^k},\dots,t^{\delta_sb_i^k}$.
Moreover, the sets $\{
t^{(\delta_1)_{b_i^k,g}b_i^k},\dots,t^{(\delta_s)_{b_i^k,g}b_i^k}\}$
for distinct values of $k$ are linearly independent by Lemma
\ref{L1}. Hence if $N\in \Z_+$ is sufficiently large and $g\prec
b_1^N,\dots,b_{2m}^N$, we obtain that $\{
t^{\alpha_1}D(g),\dots,t^{\alpha_{r_0}}D(g)\}$ is linearly
independent in $A(g)$ if and only if $\{
t^{\delta_1}D(g),\dots,t^{\delta_{r_0}}D(g)\}$ is. But $\max\{
\ell(\delta_1),\dots,\ell(\delta_{r_0})\} < \ell$, so the lemma
follows by induction. \pfend

\begin{corollary}
\label{key'}
Let $t^{\alpha_1 f_1},\dots,t^{\alpha_r f_r}\in Gen$ be distinct with
 $f_1,\dots,f_r\in G_0$ and $\alpha_1,\dots,\alpha_r\in F(\bar y)$. Then there is $g\in \Lambda$ such that $$\{ t^{\alpha_1 f_1}D(g),\dots,t^{\alpha_r f_r}D(g)\}\;\;\hbox{ is linearly independent in } \;A(g).$$ \end{corollary}

{\it Proof}: Let $f_s=v_s(\bar b)h_s$ and any non-empty $\alpha_s=:y_{i_s}^{-1}\beta_s$ ($s=1,\dots,r$). Let $N\in \Z_+$ be sufficiently large with $N>2(\ell(v_1)+\dots+\ell(v_r))$. If $g\prec b_1^{2N},\dots,b_{2m}^{2N}$ and $s\in \{ 1,\dots,r\}$ with $\alpha_s$ non-empty, then $t^{\alpha_sv_sh_s}D(g)=t^{(\alpha_s)_{f_s,g}v_sh_s}D(g)$ contains a term with the conjugator ending $b_{i_s}^Nv_sh_s$; this portion, in reduced form, begins with $b_{i_s}^{k_s}$ for some $k_s\geq N/2>\ell(v_s)$. By Lemma \ref{L1} and the choice of $N$, we deduce that all $h_s$ are equal as are all $v_s$; moreover, $\{ t^{\alpha_1 f_1}D(g),\dots,t^{\alpha_r f_r}D(g)\}$ is a linearly independent set in $A(g)$ if and only if $\{ t^{\alpha_1}D(g),\dots,t^{\alpha_r}D(g)\}$ is. The corollary now follows from Lemma \ref{key}. \pfend

We now use Nielsen's method (see, {\it e.g.}, \cite{LS}, Chapter 1 Section 2) to lift Corollary \ref{key'} to the non-abelian case.

\begin{corollary} \label{L6} Distinct elements $t^{\alpha_1f_1},\dots,t^{\alpha_rf_r}\in Gen$ generate a free subgroup. 
\end{corollary}

{\it Proof}: It is enough to prove  that there is $g\in \Lambda$ such that $u_1=t^{\alpha_1f_1}D(g)$, $\dots$,  $u_r=t^{\alpha_rf_r}D(g)$ are free generators of the free subgroup $\langle u_1,\dots,u_r\rangle$ of $T/D(g)$. Indeed, by Corollary \ref{key'}, there is $g\in \Lambda$ with
$\{ u_1,\dots,u_r\}$ linearly independent in $A(g)$. Let $U$ be the subgroup of $T$ generated by
$t^{\alpha_1f_1},\dots,t^{\alpha_rf_r}$. By Lemma~\ref{free1} (iv),
$D(g)U/D(g)$ is a subgroup of the free group $T/D(g)$ of rank {\it at
most} $r$. But $u_1,\dots,u_r$ are linearly independent in $A(g)$ and so in $U/\ga_2(U)$. So the free
group $D(g)U/D(g)$ has rank $r$. Now
$u_1,\dots,u_r$ generate the subgroup $D(g)U/D(g)$ of the free group of rank $r$ in $T/D(g)$. By
\cite{LS}, Chapter 1, Proposition 2.7,  we get that $\{u_1,\dots,u_r\}$ is a free
generating set for $D(g)U/D(g)$. \pfend

\begin{proposition}
\label{prop3}
The group $G$ is fintely presented.
\end{proposition}

{\it Proof}: We claim that $G$ has presentation (\ref{d5}), (\ref{d1}) -- (\ref{d4}) in generators $t$, $x_1,\dots, x_m$, $b_1,\dots,b_m$, $y_1,\dots,y_{2m}$. By Lemma~\ref{RelG}, we can also use relations (\ref{d6}), (\ref{top}), and (\ref{comy}) that follow from (\ref{d5}), (\ref{d1}) -- (\ref{d4}). Assume $w(t,\bar x,\bar b,\bar y)=1$ in $G$. Collecting $x_i$, $b_j$, and $y_k$ to the left and using relations (\ref{d5}), (\ref{d1}) and (\ref{d2}), we can rewrite $w$ in the form $w=w_1(\bar x)w_2(\bar b)w_3(\bar y)u$ where $u\in T$ and $w_1$, $w_2$, $w_3$ are reduced words in indicated generators. By (\ref{GT}), the words $w_1$, $w_2$, $w_3$ are empty. Thus $w\in T$. 

Using Lemma~\ref{L5}, we write $w$ in generators $t^{\alpha_1f_1},\dots, t^{\alpha_kf_k}\in Gen$. By Corrolary~\ref{L6}, the subgroup $\langle t^{\alpha_1f_1},\dots, t^{\alpha_kf_k}\rangle$ is free; so $w$ is the empty word.
\pfend

\medskip

Theorem~\ref{IWPO} now follows from Propositions~\ref{prop1}, \ref{prop2}, \ref{prop3}.

\bigskip

We observe that the above proof actually gives

\begin{theorem}
Let $\hat H$ be a right-orderable finitely presented group on $m$ generators and $\varepsilon$ be the natural homomorphism from the free group $F_m$ onto $\hat H$. Then the semidirect product $\langle t^{\hat h}\mid \hat h\in \hat H\rangle\rtimes F$ defined by automorphisms $(t^{\hat h})^f=t^{\hat h\varepsilon(f)}$ is embeddable in a (two-sided) orderable finitely presented group.
\end{theorem}

Authors' addresses:
\bigskip

V. V. Bludov:

Chair of Mathematics,

Baikal National University of Economics and Law,

Irkutsk  664011,

Russia
\medskip

vasily-bludov@yandex.ru
\bigskip

A. M. W. Glass:
\medskip

Queens' College,

Cambridge CB3 9ET,

England
\medskip

and
\medskip

Department of Pure Mathematics and Mathematical Statistics,

Centre for Mathematical Sciences,

Wilberforce Rd.,

Cambridge CB3 0WB,

England
\medskip

amwg@dpmms.cam.ac.uk
\bigskip

\end{document}